\documentclass[10pt]{article}

\usepackage{graphicx} 
\usepackage{amsmath} 
\usepackage{amssymb} 
\usepackage{amsthm}
\usepackage{yfonts}
\usepackage{mathrsfs}
\usepackage{xypic}
\title{Injectivity radius and Cartan polyhedron for simply connected symmetric spaces} 
\author{Ling Yang\\Fudan University}
\date{2006/5/22}%

\begin{document}
\def\e{\mathbf{e}}
\def\r{\mathbf{r}}
\def\om{\omega}
\def\Om{\Omega}
\def\td{\tilde}
\def\w{\wedge}
\def\c{\cdot}
\def\n{\nabla}
\def\y{\mathbf{y}}
\def\a{\mathbf{a}}
\def\b{\mathbf{b}}
\def\Y{\mathbf{Y}}
\def\p{\partial}
\def\f{\frac}
\def\si{\sigma}
\def\mc{\mathcal}
\def\C{\Bbb{C}}
\def\ra{\rightarrow}
\def\lan{\langle}
\def\ran{\rangle}
\def\i{\sqrt{-1}}
\def\la{\lambda}
\def\tr{\mbox{tr}}
\def\R{\Bbb{R}}
\def\Z{\Bbb{Z}}
\def\Q{\Bbb{Q}}
\def\ol{\overline}
\def\th{\theta}
\def\td{\tilde}
\def\e{\eta}
\def\ep{\epsilon}
\def\De{\Delta}
\def\a{\alpha}
\def\be{\beta}
\def\z{\zeta}
\def\La{\Lambda}
\def\la{\lambda}
\def\g{\gamma}
\def\de{\delta}
\def\Th{\Theta}
\def\p{\partial}
\def\k{\kappa}
\def\rank{\mbox{rank}}
\def\Id{\mathbf{Id}}
\def\d{\dot}
\def\dd{\ddot}
\def\sw{\textswab}
\def\Inn{\mbox{Inn}}
\def\wtd{\widetilde}
\def\Si{\Sigma}
\def\G{\Gamma}
\def\P{\Bbb{P}}
\def\H{\Bbb{H}}
\newtheorem{Pro}{Proposition}[section]
\newtheorem{Lem}{Lemma}[section]
\newtheorem{Thm}{Theorem}[section]
\newtheorem{Cor}{Corollary}[section]
\newtheorem{Deo}{Deotation}[section]
\newtheorem{Rem}{Remark}[section]
\renewcommand{\theequation}{\arabic{section}.\arabic{equation}}

\maketitle

\begin{abstract}
We explore relationship between the cut locus of an arbitrary simply connected and compact Riemannian symmetric space
and the Cartan polyhedron of corresponding restricted root system, and compute injectivity radius and diameter for
every type of  irreducible ones.
\newline

\newline

\textbf{2000 MR Subject Classification: 53C35}
\end{abstract}
\section*{0\ \ \  Introduction}

Let $(M,g)$ be an $n$-dimensional Riemannian manifold, it is well known that for each $p\in M$ the exponential map is injective on a sufficiently
small ball in $M_p$, then there is a natural question to be taken up: how to determine the maximum radius of such a
ball (i.e., \textbf{injectivity radius} of $M$, denoted by $i(M)$)?
Meanwhile, if we assume $M$ to be compact, then the length of an arbitrary minimal
geodesic in $M$ have a least upper bound (i.e., \textbf{diameter} of $M$, denoted by $d(M)$).
Injectivity radius and diameter have a close relationship
with curvature of $M$, which could be easily seen from Bonnet-Myers Theorem and Klingenberg Theorem: the former gives
a upper bound of $d(M)$
when the Ricci curvature of $M$  has a positive lower bound; the latter tells us $i(M)$ is no less than
$\pi$ when $M$ is simply connected, $n\geq 3$ and $1/4<K\leq 1$. Cheeger, Toponogov, Berger, Grove and Shiohama
have made a contribution to this topic (see \cite{CE} Ch. 5-6 and \cite{GS}).

The purpose of this paper is to determine injectivity radius and diameter for an arbitrary of simply connected and compact Riemannian symmetric space
explicitly. The author hopes the results be beneficial to doing further research for geometric properties on symmetric spaces of compact type.

Our computation bases on the work of
Richard Crittenden, who discussed conjugate points and cut points in symmetric spaces in \cite{Cr}. In the paper,
he claimed that the conjugate locus is determined by the diagram of a single Cartan subalgebra and the isotropy
group, and proved that the cut locus of $p$ coincides with the first conjugate locus of $p$ for every $p\in M$ using
algebraic method (Cheeger proved the same conclusion using geometrical method, see \cite{Ch}). But he didn't
describe the first conjugate locus precisely. In Section 1, we explore relationship between the first conjugate locus and the Cartan polyhedron
of corresponding restricted root system after summarizing the results due to Richard Crittenden,
then we obtain the main theorem (i.e., Theorem 1.3) about the cut locus;
our denotation is mainly from \cite{Hel} and \cite{AB}. (Theorem 1.1 is same as \cite{Hel} p.294 Prop. 3.1,
but our method is different.)

Theorem 1.3 and the definition of $i(M)$ and $d(M)$ tell us, both geometrical quantities can be determined from researching into the properties of Cartan polyhedron,
and they depend on the type of the restricted root system $\Si$, the type of the corresponding orthogonal
involutive Lie algebra, and the metric on $M$. So we divide the process of computing $i(M)$ and $d(M)$ into three
steps. Section 2 is an independent section, in which the subject we deal with is an arbitrary abstract irreducible
root system $\Phi$; we define two new variables, i.e., $i(\Phi)$ and $d(\Phi)$,
which only depend on $\Phi$, show that $i(\Phi)$ is equal to the reciprocal
of the length of the highest root, and explore the relationship between $d(\Phi)$ and the length of the highest root
for every type of irreducible root system. Section 3 is the complement of Section 2, in which the subject we discuss
is the restricted root system $\Si$, and the inner product on it is induced from the Killing form on the orthogonal
involutive Lie algebra associated with $M$; we compute the square of the length of the highest root of $\Si$,
which is depend on the type of the orthogonal involutive Lie algebra but has no direct relationship with the type
of $\Si$; the computation is on the basis of the \textbf{Satake diagram} given by Araki in \cite{Ar}. In Section 4,
the subject we research into is an arbitrary simply connected, compact and irreducible Riemannian
symmetric space; at the beginning of the section, we define a parameter $\ep>0$ which only depend on the metric of
$M$ and show the relationship between $\ep$ and the Ricci curvature of $M$; then we claim $i(M)=\pi\k^{1/2}$, where
$\k$ denotes the maximum of the sectional curvatures of $M$, and list $i(M)$ and $d(M)$ for every type
of them when $\ep=1$, $Ric=1/2$ in Table 4.1 and Table 4.2 on the basis
of what we have done in Section 2-3; at last we spend a little effort to discuss the reducible cases.

The author wishes to express his sincere gratitude to his supervisor, Professor Xin Yuan-long, for his inspiring suggestions, as well as to
Doctor Liu Xu-Sheng for providing some references.

\section{Conjugate locus and cut locus of an arbitrary compact and simply connected Riemannian symmetric space}
\setcounter{equation}{0}

Let $(M,g)$ be an arbitrary  Riemannian locally symmetric space with non-negative sectional curvature, i.e., $\n R=0$ and $K\geq 0$, where
$\n$ is Levi-Civita connection corresponding to $g$, and $R$ is corresponding curvature tensor field on $M$($R(X,Y)=-[\n_X,\n_Y]+\n_{[X,Y]}$).
For arbitrary $X\in T_o M$,
let $\g:(-\infty,\infty)\ra M$ be a geodesic satisfying $\d{\g}(0)=X$ (i.e., $\g(t)=\exp_o(tX)$).
A vector field $U$ along $\g$ is called a \textbf{Jacobi field} if it satisfies the Jacobi equation:
\begin{eqnarray}
\dd{U}+R_{\d{\g}U}\d{\g}=0.
\end{eqnarray}

Define a self-adjoint map $T_X: (T_o M,\lan,\ran)\ra (T_o M,\lan,\ran)$
$Y\mapsto R_{X,Y}X$, where $\lan Y,Z\ran=g(Y,Z)$ for every $Y,Z\in T_o M$; denote by $\la_1,\cdots,\la_m\geq 0$
the eigenvalues of $T_X$, by $(T_X)_i$ the eigenspace with respect to $\la_i$, then
\begin{eqnarray}
T_o M=\bigoplus_{i=1}^m (T_X)_i.
\end{eqnarray}
For arbitrary $Y_i\in (T_X)_i$, let
$Y_i(t)$ be the vector field obtained by parallel translation of $Y_i$ along $\g$, then the Jacobi field satisfying $U_i(0)=0$ and $\d{U_i}(0)=Y_i$ is
\begin{eqnarray}
U_i(t)=\left\{\begin{array}{ll}t Y_i(t) & \mbox{when }\la=0;\\ \f{1}{\sqrt{\la}}\sin(\sqrt{\la}t)Y_i(t) & \mbox{when }\la>0.\end{array}\right.
\end{eqnarray}
(cf. \cite{Xin} p.195). And moreover, the Jacobi field satisfying $U(0)=0$ and $\d{U}(0)=Y=\sum_{i=1}^m Y_i$, where $Y_i\in (T_X)_i$,
is $U=\sum_{i=1}^m U_i$.

$X$ is called a conjugate point in $T_o M$, if and only if there exists a nonzero Jacobi field $U$ along $\g(t)=\exp_o(tX)$, such that $U(0)=U(1)=0$.
By (1.2) and (1.3), we immediately obtain the following Proposition:

\begin{Pro}
Let $(M,g)$ be a Riemannian locally symmetric space with non-negative sectional curvature, fix $o\in M$ and $X\in T_o M$; denote $T_X(Y)=R_{X,Y}X$,
let $\la_1,\cdots,\la_m$ be the eigenvalues of $T_X$; then $X$ is a conjugate point in $T_o M$
if and only if there exists at least one positive eigenvalue $\la_i$ of $T_X$ such that $\sqrt{\la_i}\in \pi\Z$.
\end{Pro}

Denote the first conjugate locus of $o$ in $T_o M$ by $K(o)$. By the definition of $K(o)$, $X\in K(o)$ if and only if
$X$ is a conjugate point and $tX$ isn't conjugate point for every $t\in (0,1)$. Notice $T_{tX}=t^2T_X$; then applying Proposition 1.1
we have:

\begin{Pro}
The assumption and denotation is similar to Proposition 1.1, then $X\in K(o)$ if and only if $\max_{1\leq i\leq m}
\sqrt{\la_i}=\pi$.
\end{Pro}

Before applying Proposition 1.1 and 1.2 to compact symmetric spaces, we recall several basic concepts about restricted root systems.

Let $\sw{u}$ be a compact semisimple Lie algebra and $\th$ an involutive automorphism of $\sw{u}$, then $\th$ extends
uniquely to a complex involutive automorphism of $\sw{g}$, the complexification of $\sw{u}$. We have then the direct decompositions
\begin{eqnarray}
\sw{u}=\sw{k}_0\oplus \sw{p}_*;\qquad \mbox{where }\sw{k}_0=\{X\in \sw{u}:\th(X)=X\},\sw{p}_*=\{X\in \sw{u}:\th(X)=-X\}.
\end{eqnarray}
Let $\lan,\ran$ be an inner product on $\sw{p}_*$ invariant under $Ad\ \sw{k}_0$, then $(\sw{u},\th,\lan,\ran)$
is an \textbf{orthogonal involutive Lie algebra}; without loss of generality we can assume it is \textbf{reduced} (cf. \cite{AB} p. 20-21).
Let $M=U/K$ with $U$-invariant metric $g$ is a compact
Riemannian symmetric space which associates with $(\sw{u},\th,\lan,\ran)$, then there is a natural correspondence
between $(T_o M,g)$ and $(\sw{p}_*,\lan,\ran)$, where $o=eK$; in the following text we identify $T_o M$ and $\sw{p}_*$.

Let $\sw{h}_{\sw{p}_*}$ denote an arbitrary maximal abelian subspace of $\sw{p}_*$, $\sw{h}_{\sw{k}_0}$ be an abelian subalgebra
of $\sw{k}_0$ such that $\sw{h}_{\sw{k}_0}\oplus \sw{h}_{\sw{p}_*}$ is a maximal abelian subalgebra of $\sw{u}$,
and $\sw{h}$ denote the subalgebra of $\sw{g}$ generated by $\sw{h}_{\sw{k}_0}\oplus \sw{h}_{\sw{p}_*}$. Denote
$\sw{p}_0=\sqrt{-1}\sw{p}_*$, $\sw{p}=\sw{p}_*\otimes \C$, $\sw{k}=\sw{k}_0\otimes \C$, $\sw{h}_{\sw{p}_0}=\sqrt{-1}\sw{h}_{\sw{p}_*}$,
$\sw{h}_{\sw{p}}=\sw{h}_{\sw{p}_*}\otimes \C$,
then the Killing form $(,)=B(,)$ is positive on $\sqrt{-1}\sw{h}_{\sw{k}_0}\oplus \sw{h}_{\sw{p}_0}$;
let $\De$ be the root system of $\sw{g}$ with respect to $\sw{h}$, then $\sqrt{-1}\sw{h}_{\sw{k}_0}\oplus \sw{h}_{\sw{p}_0}$
is the real linear space generated by $\De$, which is denoted by $\sw{h}_\R$. Denote by $\De^+$ the subset of $\De$ formed by the positive roots with
respect to a lexicographic ordering of $\De$; for every $\a\in \De$, denote by $\a^\th=\th(\a)$, by $\bar{\a}=1/2(\a-\a^\th)$
the orthogonal projection of $\a$ into $\sw{p}_0$. Denote by $\De_0=\{\a\in \De:\bar{\a}=0\}$, $\De_{\sw{p}}=\{\a\in \De:\bar{\a}\neq 0\}$, $P_+=\De^+\cap
\De_{\sw{p}}$; by $\Si=\{\bar{\a}:\a\in \De_{\sw{p}}\}$ the \textbf{restricted root system}. $\Si$
has a compatible ordering with $\De$, and $\Si^+=\{\bar{\a}:\a\in P_+\}$. Denote by
\begin{eqnarray}
&&\sw{g}_\g=\{x\in \sw{g}:[H,x]=(H,\g)x,H\in \sw{p}\}\qquad \g\in \Si,\\
&&\sw{k}_\g=(\sw{g}_\g\oplus \sw{g}_{-\g})\cap \sw{k}, \sw{p}_\g=(\sw{g}_\g\oplus \sw{g}_{-\g})\cap \sw{p}\qquad \g\in \Si^+,
\end{eqnarray}
and by $m_\g=\dim_\C \sw{g}_\g$ the \textbf{multiplicity} of $\g$, then
\begin{eqnarray}
\sw{p}=\sw{h}_{\sw{p}}\oplus \Big(\bigoplus_{\g\in \Si^+}\sw{p}_\g\Big)
\end{eqnarray}
and
\begin{eqnarray}
m_\g=\big|\{\a\in \De_{\sw{p}}:\bar{\a}=\g\}\big|,\ \dim \sw{k}_\g=\dim \sw{p}_\g=m_\g
\end{eqnarray}
(cf. \cite{Hel} p. 283-293).

It is well known that
\begin{eqnarray}
R_{X,Y}Z=ad[X,Y]Z\qquad X,Y,Z\in \sw{p}_*
\end{eqnarray}
(cf. \cite{Ko2} p. 231, in which the curvature tensor is defined by $R(X,Y)=\n_X\n_Y-\n_Y\n_X-\n_{[X,Y]}$ ); i.e., $T_X(Y)=R_{X,Y}X=-(ad X)^2 Y$, $T_X=-(ad X)^2$.

For every $X\in \sw{p}_*$, there exists $k\in K$ and $H\in \sw{h}_{\sw{p}_*}$, such that $X=Ad(k)H$ (cf. \cite{AB} p. 31).
For arbitrary $u\in \sw{p}_\g$, $-(ad\ H)^2 u=\big(ad\ (-\sqrt{-1}H)\big)^2 u=(-\sqrt{-1}H,\g)^2 u$;
$-\sqrt{-1}H,\g\in \sw{h}_{\sw{p}_0}$ yields $(-\sqrt{-1}H,\g)^2\geq 0$; then by (1.5) and (1.7), the eigenvalues of $T_H=-(ad\ H)^2$ include
\begin{eqnarray}
0,\ (-\sqrt{-1}H,\g)^2\ (\g\in \Si^+).
\end{eqnarray}
Since $X=Ad(k)H$, $ad\ X=Ad(k)\circ ad\ H\circ Ad(k)^{-1}$ and moreover
$T_X=Ad(k)\circ T_H\circ Ad(k)^{-1}$; which yields the eigenvalues of
$T_X$ coincide with the eigenvalues of $T_H$. Applying Proposition 1.1, we have

\begin{Thm}
Let $M=U/K$ be a compact Riemannian symmetric space such that $U$ is a semi-simple and compact Lie group,
and the denotation of $\sw{p}_*,\sw{k}_0,\sw{h}_{p_*},\Si$ is similar to above, then
for every $X=Ad(k)H\in \sw{p}_*$, where $k\in K,H\in \sw{h}_{p_*}$, $X$ is a conjugate point in $T_o M$
if and only if there exists at least one $\g\in \Si$, such that
$(H,\g)\in \pi\sqrt{-1}(\Z-0)$.
\end{Thm}

Now we denote by $C$ the Weyl chamber with respect the ordering of $\Si$, i.e., $C=\{x\in \sw{h}_{\sw{p}_0}:
(x,\g)>0\mbox{ for every }\g\in \Si^+\}$, by $\Pi$ the set of simple roots.
Recall that the planes
$(x,\g)\in \Z$($\g\in \Si$) in $\sw{h}_{\sw{p}_0}$ constitute the \textbf{diagram} $D(\Si)$ of $\Si$, and the closure of a connected
component of $\sw{h}_{\sw{p}_0}-D(\Si)$ will be called a \textbf{Cartan polyhedron}. Especially,
let $A$ be the set of maximal roots, then the inequalities $(x,\g)\geq 0(\g\in \Pi), (x,\be)\leq 1(\be\in B)$ define
a Cartan polyhedron, which is denoted by $\triangle$ (See \cite{AB} p. 10). Obviously $\triangle\subset \ol{C}$, where $\ol{C}$
denotes the closure of $C$ in $\sw{h}_{\sw{p}_0}$. Since
Weyl group $W$ permutes Weyl chamber in a simply transitive manner  and every element of Weyl group
can be extended to $Ad_{\sw{u}}(\sw{k}_0)$ (See \cite{Hel} p. 288-290), for every $X\in \sw{p}_*$, there exists
$k\in K$ and $H\in \sqrt{-1}\ \ol{C}$ such that $X=Ad(k)H$. By Proposition 1.2 and (1.10), $X\in K(o)$ if and only if
\begin{eqnarray}
\pi=\max_{\g\in \Si^+}\big|(-\sqrt{-1}H,\g)\big|=\max_{\be\in B}\big|(-\sqrt{-1}H,\be)\big|;\qquad (\mbox{ since }-\sqrt{-1}H\in \ol{C})
\end{eqnarray}
i.e., $H\in \pi\sqrt{-1}\triangle$ and $(-\sqrt{-1}H,\be)=\pi$ for some $\be\in B$.

As a matter of convenience, we bring in new denotation:

\begin{Deo}
Denote $\triangle'=\{x\in \triangle:(x,\be)=1\mbox{ for some }\be\in B\}$, i.e., $\triangle'$ is the union of the
facets of $\triangle$ which don't contain $0$.
\end{Deo}

Then we have:

\begin{Thm}
The assumption is same to Theorem 1.1, then $K(o)=Ad(K)(\pi\sqrt{-1}\triangle')$.
\end{Thm}

On cut locus,
in 1962, Richard Crittenden proved the following proposition (See \cite{Cr}):

\begin{Lem}
Let $M$ be a simply connected complete symmetric space, for every $p\in M$, the cut locus of
$p$ coincides with the first conjugate locus of $p$.
\end{Lem}

For every $p\in M$, denote by $C(p)$ the cut locus of $p$ in $T_p M$.
Let $F:M\ra M$ be an isometry, then for any $p\in M$ and $X\in T_p M$, $d\big(p,\exp_p(X)\big)=|X|$ if and only if
$d\big(F(p),\exp_{F(p)}((dF)_p X)\big)=\big|(dF)_pX\big|$, which yields $C(q)=(dF)_p C(p)$. Then by Theorem 1.2 and
Lemma 1.1, we have

\begin{Thm}
Let $M=U/K$ be a simply connected and compact Riemannian symmetric space such that $U$ is a semi-simple and compact Lie group,
and the denotation of $\sw{p}_*,\sw{k}_0,\sw{h}_{p_*},\Si,\triangle'$ is similar to above, then the cut locus of $o$
in $T_o M=\sw{p}_*$ is $C(o)=Ad(K)(\pi\sqrt{-1}\triangle')$; and for any $p=aK\in M$, $C(p)=(dL_a)_o C(o)$, where $a\in U$ and
$L_a$ is an isometry of $M$ satisfying $L_a(bK)=abK$.
\end{Thm}

\section{Some computation on Cartan polyhedron}
\setcounter{equation}{0}

In the section, we assume $\Phi\subset V$ be a irreducible abstract root system with an ordering, where $V$ is an
$l$-dimensional real vector space with inner product $(,)$; denote by $\Phi^+$
and $\Pi=\{\a_1,\cdots,\a_l\}$ respectively the positive root system and the simple root system, and by $\psi$ the highest root;
let $d_1,\cdots,d_l\in \Z^+$ such that $\psi=\sum_{i=1}^l d_i \a_i$(cf. \cite{AB} p. 9-10). The definition of $\triangle$ and $\triangle'$
is similar to Section 1 (notice that in this case $A=\{\psi\}$);
since $\Phi$ is irreducible, $\triangle$ is a simplex, whose vertices are $0,e_1,\cdots,e_l$, which satisfy
\begin{eqnarray}
(e_j,\a_i)=\f{1}{d_j}\de_{ij}.
\end{eqnarray}
Define
\begin{eqnarray}
i(\Phi)=\min_{x\in \triangle'}(x,x)^{1/2},\ d(\Phi)=\max_{x\in \triangle'}(x,x)^{1/2};
\end{eqnarray}
and in the following we will compute $i(\Phi)$ and $d(\Phi)$.

For every $x\in \triangle'$,
$1=(x,\psi)\leq (x,x)^{1/2}(\psi,\psi)^{1/2}$, and
the equal sign holds if and only if $x=\psi/(\psi,\psi)\in \triangle'$ (since $(\psi,\a_i)\geq 0$ for every $\a_i\in \Pi$),
thus
\begin{eqnarray}
i(\Phi)=(\psi,\psi)^{-1/2}.
\end{eqnarray}

$x\mapsto (x,x)^{1/2}$ is a function on $\triangle'$; since for any $t\in (0,1)$,
\begin{eqnarray}
(tx_1+(1-t)x_2,tx_1+(1-t)x_2)^{1/2}&=&\big(t^2(x_1,x_1)+(1-t)^2(x_2,x_2)+2t(1-t)(x_1,x_2)\big)^{1/2}\nonumber\\
&\leq& t(x_1,x_1)^{1/2}+(1-t)(x_2,x_2)^{1/2};
\end{eqnarray}
the function reach its maximum at the vertices of $\triangle'$, including $e_1,\cdots,e_l$. Denote $\Om_{ij}=(\a_i,\a_j)$ and $e_j=\a_k A_j^k$;
(2.1) yields
$$\f{1}{d_j}\de_{ij}=(e_j,\a_i)=(\a_k A_j^k,\a_i)=\Om_{ik}A_j^k;$$
so $A_j^k=1/d_j(\Om^{-1})_{ki}\de_{ij}=1/d_j(\Om^{-1})_{kj}$ and
\begin{eqnarray}
&&(e_i,e_j)=(\a_k A_i^k,e_j)=\f{1}{d_j}\de_{jk}A_i^k=\f{1}{d_j}\de_{jk}\f{1}{d_i}(\Om^{-1})_{ki}=\f{1}{d_jd_i}(\Om^{-1})_{ji};\nonumber\\
&&\mbox{especially }(e_j,e_j)=\f{1}{d_j^2}(\Om^{-1})_{jj}.
\end{eqnarray}
Thus
\begin{eqnarray}
d(\Phi)=\max_{1\leq j\leq l}(e_j,e_j)^{1/2}=\max_{1\leq j\leq l}\f{1}{d_j}{(\Om^{-1})_{jj}}^{1/2}.
\end{eqnarray}

A root system $\Phi$ is said to be \textbf{reduced}, if and only if for every $\a,\be\in \Phi$ which are proportional,
we have $\a=\pm \be$. It is well known that the root systems $\sw{a}_l(l\geq 1)$, $\sw{b}_l(l\geq 2)$, $\sw{c}_l(l\geq 3)$,
$\sw{d}_l(l\geq 4)$, $\sw{e}_6$, $\sw{e}_7$, $\sw{e}_8$, $\sw{f}_4$, $\sw{g}_2$ exhaust all irreducible reduced root systems,
and every irreducible reduced root system associates with a unique Dynkin diagram. Otherwise, every irreducible nonreduced
root system is isomorphic to $(\sw{bc})_l(l\geq 1)$; the set of \textbf{indivisible roots} ($\a\in \Phi$ is called indivisible if and only
if $1/2\a\notin \Phi$) in $(\sw{bc})_l$ is isomorphic to $\sw{b}_l$ (cf. \cite{Hel} p. 474-475). In the following we
give the detail of computing $d(\Phi)$
for every type of root systems.
\newline

$\Phi=\sw{a}_l$: The corresponding Dynkin diagram is
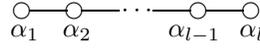
\setlength{\unitlength}{1mm}
\begin{center}
\begin{picture}(33,5)
\put(1.5,3.5){\circle{2}}
\put(2.5,3.5){\line(1,0){5}}
\put(8.5,3.5){\circle{2}}
\put(9.5,3.5){\line(1,0){5}}
\put(14.8,2.6){$\cdots$}
\put(19,3.5){\line(1,0){5}}
\put(25,3.5){\circle{2}}
\put(26,3.5){\line(1,0){5}}
\put(32,3.5){\circle{2}}
\put(0,0){$\a_1$}
\put(7,0){$\a_2$}
\put(21,0){$\a_{l-1}$}
\put(30.5,0){$\a_l$}
\end{picture}
\end{center}
then $\psi=\sum_{i=1}^l \a_i$, $d_i=1$ for every $i$.
Denote $\a_1=x_1-x_2,\cdots, \a_l=x_l-x_{l+1}$, then $\psi=x_1-x_{l+1}$ and therefore
$(x_i,x_j)=1/2(\psi,\psi)\de_{ij}$;
by (2.1) and (2.6), we obtain
\begin{eqnarray}
&&e_j=\f{2}{(\psi,\psi)(l+1)}\big((l+1-j)\sum_{k=1}^j x_k-j\sum_{k=j+1}^{l+1}x_k\big)\qquad 1\leq j\leq l;\\
&&d(\sw{a}_l)=\max_{1\leq j\leq l} (e_j,e_j)^{1/2}=\left\{\begin{array}{ll}\f{\sqrt{2}}{2}(\psi,\psi)^{-1/2}(l+1)^{1/2} & l\mbox{ is odd;}\\\f{\sqrt{2}}{2}(\psi,\psi)^{-1/2}\big(l(l+2)\big)^{1/2}(l+1)^{-1/2} & l\mbox{ is even.}\end{array}\right.
\end{eqnarray}

$\Phi=\sw{b}_l$: The corresponding Dynkin diagram is
\setlength{\unitlength}{1mm}
\begin{center}
\begin{picture}(34,5)
\put(1.5,3.5){\circle{2}}
\put(2.5,3.5){\line(1,0){5}}
\put(8.5,3.5){\circle{2}}
\put(9.5,3.5){\line(1,0){5}}
\put(14.8,2.6){$\cdots$}
\put(19,3.5){\line(1,0){5}}
\put(25,3.5){\circle{2}}
\put(25.9,2.6){$\Longrightarrow$}
\put(32.2,3.5){\circle{2}}
\put(0,0){$\a_1$}
\put(7,0){$\a_2$}
\put(21,0){$\a_{l-1}$}
\put(30.5,0){$\a_l$}
\end{picture}
\end{center}
then $\psi=\a_1+2\sum_{i=2}^l \a_i$, $d_1=1$ and $d_i=2$ for every $2\leq i\leq l$. Denote $\a_1=x_1-x_2,\cdots,\a_{l-1}
=x_{l-1}-x_l, \a_l=x_l$, then $\psi=x_1+x_2$ and therefore
$(x_i,x_j)=1/2(\psi,\psi)\de_{ij}$;
by (2.1) and (2.6), we obtain
\begin{eqnarray}
&&e_1=\f{2}{(\psi,\psi)}x_1,\ e_j=\f{1}{(\psi,\psi)}\sum_{k=1}^j x_k\ (2\leq j\leq l);\\
&&d(\sw{b}_l)=\max_{1\leq j\leq l} (e_j,e_j)^{1/2}=\left\{\begin{array}{ll} \sqrt{2}(\psi,\psi)^{-1/2} & l\leq 3;\\\f{\sqrt{2}}{2}(\psi,\psi)^{-1/2}l^{1/2} & l\geq 4.\end{array}\right.
\end{eqnarray}

$\Phi=\sw{c}_l$: The corresponding Dynkin diagram is
\setlength{\unitlength}{1mm}
\begin{center}
\begin{picture}(34,5)
\put(1.5,3.5){\circle{2}}
\put(2.5,3.5){\line(1,0){5}}
\put(8.5,3.5){\circle{2}}
\put(9.5,3.5){\line(1,0){5}}
\put(14.8,2.6){$\cdots$}
\put(19,3.5){\line(1,0){5}}
\put(25,3.5){\circle{2}}
\put(25.9,2.6){$\Longleftarrow$}
\put(32.2,3.5){\circle{2}}
\put(0,0){$\a_1$}
\put(7,0){$\a_2$}
\put(21,0){$\a_{l-1}$}
\put(30.5,0){$\a_l$}
\end{picture}
\end{center}
then $\psi=2\sum_{i=1}^{l-1} \a_i+\a_l$, $d_l=1$ and $d_i=2$ for every $1\leq i\leq l-1$. Denote $\a_1=x_1-x_2,\cdots,\a_{l-1}
=x_{l-1}-x_l, \a_l=2x_l$, then $\psi=2x_1$ and therefore
$(x_i,x_j)=1/4(\psi,\psi)\de_{ij}$;
by (2.1) and (2.6), we obtain
\begin{eqnarray}
&&e_j=\f{2}{(\psi,\psi)}\sum_{k=1}^j x_k\qquad 1\leq j\leq l;\\
&&d(\sw{c}_l)=\max_{1\leq j\leq l} (e_j,e_j)^{1/2}=(\psi,\psi)^{-1/2}l^{1/2}.
\end{eqnarray}

$\Phi=\sw{d}_l$: The corresponding Dynkin diagram is
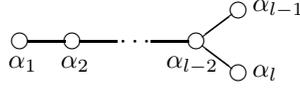
\setlength{\unitlength}{1mm}
\begin{center}
\begin{picture}(33,11.5)
\put(1.5,6){\circle{2}}
\put(2.5,6){\line(1,0){5}}
\put(8.5,6){\circle{2}}
\put(9.5,6){\line(1,0){5}}
\put(14.8,5.1){$\cdots$}
\put(19,6){\line(1,0){5}}
\put(25,6){\circle{2}}
\put(0,2.5){$\a_1$}
\put(7,2.5){$\a_2$}
\put(21,2.5){$\a_{l-2}$}
\put(25.8,6.6){\line(4,3){4}}
\put(30.6,10.2){\circle{2}}
\put(25.8,5.4){\line(4,-3){4}}
\put(30.6,1.8){\circle{2}}
\put(32.5,10){$\a_{l-1}$}
\put(32.5,1.3){$\a_l$}
\end{picture}
\end{center}
then $\psi=\a_1+2\sum_{i=2}^{l-2} \a_i+\a_{l-1}+\a_l$, $d_1=d_{l-1}=d_l=1$ and $d_i=2$ for every $2\leq i\leq l-2$. Denote $\a_1=x_1-x_2,\cdots,\a_{l-1}
=x_{l-1}-x_l, \a_l=x_{l-1}+x_l$, then $\psi=x_1+x_2$ and therefore
$(x_i,x_j)=1/2(\psi,\psi)\de_{ij}$;
by (2.1) and (2.6), we obtain
\begin{eqnarray}
&&e_1=\f{2}{(\psi,\psi)}x_1,\ e_j=\f{1}{(\psi,\psi)}\sum_{k=1}^j x_k\ (2\leq j\leq l-2),\nonumber\\
&&e_{l-1}=\f{1}{(\psi,\psi)}(\sum_{k=1}^{l-1}x_k-x_l),\ e_l=\f{1}{(\psi,\psi)}\sum_{k=1}^l x_k;\\
&&d(\sw{d}_l)=\max_{1\leq j\leq l} (e_j,e_j)^{1/2}=\f{\sqrt{2}}{2}(\psi,\psi)^{-1/2}l^{1/2}.
\end{eqnarray}

$\Phi=\sw{e}_6$: The corresponding Dynkin diagram is
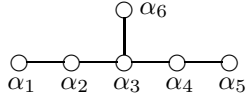
\setlength{\unitlength}{1mm}
\begin{center}
\begin{picture}(31,11.5)
\put(1.5,3.5){\circle{2}}
\put(2.5,3.5){\line(1,0){5}}
\put(8.5,3.5){\circle{2}}
\put(9.5,3.5){\line(1,0){5}}
\put(15.5,3.5){\circle{2}}
\put(16.5,3.5){\line(1,0){5}}
\put(22.5,3.5){\circle{2}}
\put(23.5,3.5){\line(1,0){5}}
\put(29.5,3.5){\circle{2}}
\put(0,0){$\a_1$}
\put(7,0){$\a_2$}
\put(14,0){$\a_3$}
\put(21,0){$\a_4$}
\put(28,0){$\a_5$}
\put(15.5,4.5){\line(0,1){5}}
\put(15.5,10.5){\circle{2}}
\put(17.5,10){$\a_6$}
\end{picture}
\end{center}
then $\psi=\a_1+2\a_2+3\a_3+2\a_4+\a_5+2\a_6$, $d_1=d_{5}=1$, $d_2=d_4=d_6=2$ and $d_3=3$.
Since all the roots have the same length,
\begin{eqnarray*}
\Om=\f{1}{2}(\psi,\psi)\left(\begin{array}{cccccc}2&-1& & & &\\-1& 2&-1& & & \\& -1& 2&-1 & & -1\\& & -1& 2& -1& \\& & & -1 & 2& \\ & & -1 & & & 2\end{array}\right);
\end{eqnarray*}
then by (2.6),
\begin{eqnarray}
d(\sw{e}_6)=\max_{1\leq j\leq 6}\f{1}{d_j}{(\Om^{-1})_{jj}}^{1/2}=\f{2\sqrt{6}}{3}(\psi,\psi)^{-1/2}.
\end{eqnarray}

$\Phi=\sw{e}_7$: The corresponding Dynkin diagram is
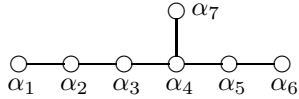
\setlength{\unitlength}{1mm}
\begin{center}
\begin{picture}(38,11.5)
\put(1.5,3.5){\circle{2}}
\put(2.5,3.5){\line(1,0){5}}
\put(8.5,3.5){\circle{2}}
\put(9.5,3.5){\line(1,0){5}}
\put(15.5,3.5){\circle{2}}
\put(16.5,3.5){\line(1,0){5}}
\put(22.5,3.5){\circle{2}}
\put(23.5,3.5){\line(1,0){5}}
\put(29.5,3.5){\circle{2}}
\put(30.5,3.5){\line(1,0){5}}
\put(36.5,3.5){\circle{2}}
\put(0,0){$\a_1$}
\put(7,0){$\a_2$}
\put(14,0){$\a_3$}
\put(21,0){$\a_4$}
\put(28,0){$\a_5$}
\put(35,0){$\a_6$}
\put(22.5,4.5){\line(0,1){5}}
\put(22.5,10.5){\circle{2}}
\put(24.5,10){$\a_7$}
\end{picture}
\end{center}
then $\psi=\a_1+2\a_2+3\a_3+4\a_4+3\a_5+2\a_6+2\a_7$, $d_1=1$, $d_2=d_6=d_7=2$, $d_3=d_5=3$ and $d_4=4$.
Since all the roots have the same length,
\begin{eqnarray*}
\Om=\f{1}{2}(\psi,\psi)\left(\begin{array}{ccccccc}2&-1& & & & &\\-1& 2&-1& & &  &\\& -1 & 2 & -1 & & &\\ & & -1& 2&-1 & & -1\\ & & & -1& 2& -1& \\& & & & -1 & 2& \\& & & -1 & & & 2\end{array}\right);
\end{eqnarray*}
then by (2.6),
\begin{eqnarray}
d(\sw{e}_7)=\max_{1\leq j\leq 7}\f{1}{d_j}{(\Om^{-1})_{jj}}^{1/2}=\sqrt{3}(\psi,\psi)^{-1/2}.
\end{eqnarray}

$\Phi=\sw{e}_8$: The corresponding Dynkin diagram is
\setlength{\unitlength}{1mm}
\begin{center}
\begin{picture}(45,11.5)
\put(1.5,3.5){\circle{2}}
\put(2.5,3.5){\line(1,0){5}}
\put(8.5,3.5){\circle{2}}
\put(9.5,3.5){\line(1,0){5}}
\put(15.5,3.5){\circle{2}}
\put(16.5,3.5){\line(1,0){5}}
\put(22.5,3.5){\circle{2}}
\put(23.5,3.5){\line(1,0){5}}
\put(29.5,3.5){\circle{2}}
\put(30.5,3.5){\line(1,0){5}}
\put(36.5,3.5){\circle{2}}
\put(37.5,3.5){\line(1,0){5}}
\put(43.5,3.5){\circle{2}}
\put(0,0){$\a_1$}
\put(7,0){$\a_2$}
\put(14,0){$\a_3$}
\put(21,0){$\a_4$}
\put(28,0){$\a_5$}
\put(35,0){$\a_6$}
\put(42,0){$\a_7$}
\put(29.5,4.5){\line(0,1){5}}
\put(29.5,10.5){\circle{2}}
\put(31.5,10){$\a_8$}
\end{picture}
\end{center}
then $\psi=2\a_1+3\a_2+4\a_3+5\a_4+6\a_5+4\a_6+2\a_7+3\a_8$, $d_1=d_7=2$, $d_2=d_8=3$, $d_3=d_6=4$, $d_4=5$ and $d_5=6$.
Since all the roots have the same length,
\begin{eqnarray*}
\Om=\f{1}{2}(\psi,\psi)\left(\begin{array}{cccccccc}2&-1& & & & & &\\-1& 2&-1& & & & &\\& -1 & 2 & -1 & & & &\\& & -1 & 2 & -1 & & &  \\& & & -1& 2&-1 & & -1\\& & & & -1& 2& -1& \\& & & & & -1 & 2& \\& & & & -1 & & & 2\end{array}\right);
\end{eqnarray*}
then by (2.6),
\begin{eqnarray}
d(\sw{e}_8)=\max_{1\leq j\leq 8}\f{1}{d_j}{(\Om^{-1})_{jj}}^{1/2}=\sqrt{2}(\psi,\psi)^{-1/2}.
\end{eqnarray}

$\Phi=\sw{f}_4$: The corresponding Dynkin diagram is
\setlength{\unitlength}{1mm}
\begin{center}
\begin{picture}(24,5)
\put(1.5,3.5){\circle{2}}
\put(2.5,3.5){\line(1,0){5}}
\put(8.5,3.5){\circle{2}}
\put(9.4,2.6){$\Longrightarrow$}
\put(15.7,3.5){\circle{2}}
\put(16.7,3.5){\line(1,0){5}}
\put(22.7,3.5){\circle{2}}
\put(0,0){$\a_1$}
\put(7,0){$\a_2$}
\put(14,0){$\a_3$}
\put(21,0){$\a_4$}
\end{picture}
\end{center}
then $\psi=2\a_1+3\a_2+4\a_3+2\a_4$, $d_1=d_4=2$, $d_2=3$ and $d_3=4$.
Since $(\psi,\psi)=(\a_1,\a_1)=(\a_2,\a_2)=2(\a_3,\a_3)=2(\a_4,\a_4)$,
\begin{eqnarray*}
\Om=\f{1}{4}(\psi,\psi)\left(\begin{array}{cccc}4 & -2 & & \\-2& 4& -2 & \\& -2& 2& -1\\& & -1 & 2\end{array}\right);
\end{eqnarray*}
then by (2.6),
\begin{eqnarray}
d(\sw{f}_4)=\max_{1\leq j\leq 4}\f{1}{d_j}{(\Om^{-1})_{jj}}^{1/2}=\sqrt{2}(\psi,\psi)^{-1/2}.
\end{eqnarray}

$\Phi=\sw{g}_2$: The corresponding Dynkin diagram is
\setlength{\unitlength}{1mm}
\begin{center}
\begin{picture}(10,5)
\put(1.5,3.5){\circle{2}}
\put(3.1,2.6){$\Rrightarrow$}
\put(7.7,3.5){\circle{2}}
\put(0,0){$\a_1$}
\put(6.2,0){$\a_2$}
\end{picture}
\end{center}
then $\psi=2\a_1+3\a_2$, $d_1=2$ and $d_2=3$. Since $(\psi,\psi)=(\a_1,\a_1)=3(\a_2,\a_2)$,
\begin{eqnarray*}
\Om=\f{1}{6}(\psi,\psi)\left(\begin{array}{cc}6 & -3\\-3 & 2\end{array}\right);
\end{eqnarray*}
and by (2.6),
\begin{eqnarray}
d(\sw{g}_2)=\max_{1\leq j\leq 2}\f{1}{d_j}{(\Om^{-1})_{jj}}^{1/2}=\f{2\sqrt{3}}{3}(\psi,\psi)^{-1/2}.
\end{eqnarray}

$\Phi=(\sw{bc})_l$: The corresponding Dynkin diagram of the set of indivisible roots in $\Phi$ is
\setlength{\unitlength}{1mm}
\begin{center}
\begin{picture}(45,5)
\put(1.5,3.5){\circle{2}}
\put(2.5,3.5){\line(1,0){5}}
\put(8.5,3.5){\circle{2}}
\put(9.5,3.5){\line(1,0){5}}
\put(14.8,2.6){$\cdots$}
\put(19,3.5){\line(1,0){5}}
\put(25,3.5){\circle{2}}
\put(25.9,2.6){$\Longrightarrow$}
\put(32.2,3.5){\circle{2}}
\put(0,0){$\a_1$}
\put(7,0){$\a_2$}
\put(21,0){$\a_{l-1}$}
\put(30.5,0){$\a_l$}
\put(40,2.3){when $l\geq 2$;}
\end{picture}
\setlength{\unitlength}{1mm}
\end{center}
\begin{center}
\begin{picture}(10,5)
\put(1.5,3.5){\circle{2}}
\put(0,0){$\a_1$}
\put(7,2.3){when $l=1$.}
\end{picture}
\end{center}
Then $\psi=2\sum_{i=1}^l \a_i$,  $d_i=2$ for every $1\leq i\leq l$. When $l\geq 2$, denote $\a_1=x_1-x_2,\cdots,\a_{l-1}
=x_{l-1}-x_l, \a_l=x_l$, then $\psi=2x_1$ and therefore
$(x_i,x_j)=1/4(\psi,\psi)\de_{ij}$;
by (2.1) and (2.6), we obtain
\begin{eqnarray}
&&e_j=\f{2}{(\psi,\psi)}\sum_{k=1}^j x_k\qquad 2\leq j\leq l;\\
&&d((\sw{bc})_l)=\max_{1\leq j\leq l} (e_j,e_j)^{\f{1}{2}}=(\psi,\psi)^{-1/2}l^{1/2}.
\end{eqnarray}
When $l=1$, $\psi=2\a_1$, $e_1=(\psi,\psi)^{-1}\psi$; so $d((\sw{bc}_1))=(e_1,e_1)^{1/2}=(\psi,\psi)^{-1/2}$; the result
coincides with (2.21).

\section{The square of the length of the highest restricted root}
\setcounter{equation}{0}

In this section, we assume $(\sw{u},\th,\lan,\ran)$ be a reduced, irreducible, semi-simple and compact orthogonal involutive Lie
algebra, $(,)=B(,)$ be the Killing form on $\sw{u}$; the denotation of
$\De,\Si,\sw{g},\sw{h},\sw{h}_\R,\sw{h}_{\sw{p}_0},\De_0,m_\g(\g\in \Si)$ is same as Section 1, and denote by $l$ and
$l_1$ respectively the rank of $\De$ and $\Si$, by $\psi\in \Si$ the highest restricted root. Then $(\sw{u},\th,\lan,\ran)$ belongs to one of the
two following types: (I) $\sw{u}$ is compact and simple, $\th$ is an involution;
(II) $\sw{u}$ is a product of two compact simple algebras exchanged by $\th$ (See \cite{AB} p. 28).
\newline

\textbf{Type I}: In the case, $\De$ and $\Si$ are both irreducible; denote by $\de$ the highest root of $\De$;
since the orderings of $\De$ and $\Si$ are compatible (i.e., $\a\geq \be$ yields $\bar{\a}\geq \bar{\be}$ for
arbitrary $\a,\be\in \De$), $\bar{\de}$ is the highest root of $\Si$; i.e., $\psi=\bar{\de}$.

Denote by $\de^\perp=\{x\in \sw{h}_{\R}:(x,\de)=0\}$, then $\De\cap \de^\perp$ is obviously a subsystem of $\De$
with an induced ordering. Denote by $\Pi=\{\a_1,\cdots,\a_l\}$ the set of simple roots in $\De$; since $(\a_i,\de)\geq 0$,
$\a=\sum_{i=1}^l a_i \a_i\in \De\cap \de^\perp$ if and only if $a_j=0$ for every $\a_j\notin \Pi\cap \de^\perp$;
which yields $\Pi\cap \de^\perp$ is the simple root system of $\De\cap \de^\perp$.

$\a_i\in \Pi\cap \de^\perp$ if and only if $\de-\a_i\notin \De\cup \{0\}$, then from the Dynkin diagram of $\De$ (which
have described in Section 2),
we can clarify $\Pi\cap \de^\perp$ and $\De\cap \de^\perp$:
\begin{eqnarray}
&&\De=\sw{a}_l:\Pi\cap \de^\perp=\{\a_i:i\neq 1,l\},\De\cap \de^\perp=\sw{a}_{l-2}\ (\emptyset \mbox{ when }l=1,2);\nonumber\\
&&\De=\sw{b}_l:\Pi\cap \de^\perp=\{\a_i:i\neq 2\},\De\cap \de^\perp=\sw{a}_{1}\oplus \sw{b}_{l-2}\ (\sw{a}_1 \mbox{ when }l=2,\sw{a}_1\oplus \sw{a}_1\mbox{ when }l=3);\nonumber\\
&&\De=\sw{c}_l:\Pi\cap \de^\perp=\{\a_i:i\neq 1\},\De\cap \de^\perp=\sw{c}_{l-1}\ (\sw{b}_2 \mbox{ when }l=3);\nonumber\\
&&\De=\sw{d}_l:\Pi\cap \de^\perp=\{\a_i:i\neq 2\},\De\cap \de^\perp=\sw{a}_1\oplus \sw{d}_{l-2}\ (\sw{a}_1\oplus \sw{a}_1\oplus \sw{a}_1 \mbox{ when }l=4,\sw{a}_1\oplus \sw{a}_3\mbox{ when }l=5);\nonumber
\end{eqnarray}
\begin{eqnarray}
&&\De=\sw{c}_l:\Pi\cap \de^\perp=\{\a_i:i\neq 1\},\De\cap \de^\perp=\sw{c}_{l-1}\ (\sw{b}_2 \mbox{ when }l=3);\nonumber\\
&&\De=\sw{d}_l:\Pi\cap \de^\perp=\{\a_i:i\neq 2\},\De\cap \de^\perp=\sw{a}_1\oplus \sw{d}_{l-2}\ (\sw{a}_1\oplus \sw{a}_1\oplus \sw{a}_1 \mbox{ when }l=4,\sw{a}_1\oplus \sw{a}_3\mbox{ when }l=5);\nonumber\\
&&\De=\sw{e}_6:\Pi\cap \de^\perp=\{\a_i:i\neq 6\},\De\cap \de^\perp=\sw{a}_5;\nonumber\\
&&\De=\sw{e}_7:\Pi\cap \de^\perp=\{\a_i:i\neq 6\},\De\cap \de^\perp=\sw{d}_6;\nonumber\\
&&\De=\sw{e}_8:\Pi\cap \de^\perp=\{\a_i:i\neq 1\},\De\cap \de^\perp=\sw{e}_7;\nonumber\\
&&\De=\sw{f}_4:\Pi\cap \de^\perp=\{\a_i:i\neq 1\},\De\cap \de^\perp=\sw{c}_3;\nonumber\\
&&\De=\sw{g}_2:\Pi\cap \de^\perp=\{\a_2\},\De\cap \de^\perp=\sw{a}_1.
\end{eqnarray}
On $\De\cap \de^\perp$, we have the following lemmas:

\begin{Lem}
$(\de,\de)=4(|\De|-|\De\cap \de^\perp|+6)^{-1}$.
\end{Lem}

\textbf{Proof}: For every $\a\in \De^+-\{\de\}$,
$\a+\de>\de$ and $\a-2\de<-\de$, which yields $p_{\a,\de}\leq 1$ and $q_{\a,\de}=0$ (for arbitrary $\a,\be\in \De$, the $\be$-string
through $\a$ is denoted by $\a-p_{\a,\be},\cdots,\a,\cdots,\a+q_{\a,\be}\be$, where $p_{\a,\be},q_{\a,\be}\in \Z^+$ satisfy $p_{\a,\be}-q_{\a,\be}
=2(\a,\be)/(\be,\be)$; cf. \cite{AB} p. 9-10);  thus
$$\f{2(\a,\de)}{(\de,\de)}=p_{\a,\de}-q_{\a,\de}=1\mbox{ or }0,\ \f{2(\a,\de)}{(\de,\de)}=1\mbox{ if and only if }(\a,\de)\neq 0\mbox{ i.e.,}\a\notin \De\cap \de^\perp.$$
By the definition of Killing form,
\begin{eqnarray*}
(\de,\de)&=&\tr(ad\ \de)^2|_{\sw{g}}=\sum_{\a\in \De}\dim \sw{g}_{\a}(\a,\de)^2=\sum_{\a\in \De}(\a,\de)^2\nonumber\\
&=&\sum_{\a\in \De-(\De\cap\de^\perp)}(\a,\de)^2=2(\de,\de)^2+\f{1}{4}(|\De|-|\De\cap\de^\perp|-2)(\de,\de)^2;
\end{eqnarray*}
that is $(\de,\de)=4(|\De|-|\De\cap \de^\perp|+6)^{-1}$. \textbf{Q.E.D}.

\begin{Lem}
$(\bar{\de},\bar{\de})=(\de,\de)$ or $1/2(\de,\de)$, and the following conditions are
equivalent:

(a) $(\bar{\de},\bar{\de})=(\de,\de)$;

(b) $\de^\th=-\de$;

(c) $\Pi_0\subset \Pi\cap \de^\perp$, where $\Pi_0=\Pi\cap \De_0$;

(d) $m_{\bar{\de}}=1$.
\end{Lem}

\textbf{Proof}: If $\de^\th=-\de$, then $\bar{\de}=\de$ and $(\bar{\de},\bar{\de})=(\de,\de)$. Otherwise
$\de^\th\neq -\de$, i.e., $\de+\de^\th
\neq 0$; Araki proved $\a+\a^\th\notin \De$ for every $\a\in \De$ in \cite{Ar}, especially $\de+\de^\th\notin \De$;
if $\de-\de^\th=2\bar{\de}\in \De\cup \{0\}$,
then $2\bar{\de}\in \Si^+$, which contradicts the assumption that $\bar{\de}$ is the highest root of $\Si$; so $(\de,\de^\th)=0$ and
$$(\bar{\de},\bar{\de})=\big(\f{\de-\de^\th}{2},\f{\de-\de^\th}{2}\big)=\f{1}{4}\big((\de,\de)+(\de^\th,\de^\th)\big)=\f{1}{2}(\de,\de).$$
So $(\bar{\de},\bar{\de})=(\de,\de)$ or $1/2(\de,\de)$ and we have proved $(a)\Longleftrightarrow (b)$.

$(b)\Longrightarrow (c)$: Suppose there exists $\a_i\in \Pi_0$ satisfying $\a_i\notin \Pi\cap \de^\perp$, then
$\bar{\a}_i=0$ and $(\de,\a_i)\neq 0$; which yields $\de-\a_i\in \De^+\cup \{0\}$ and $\th(\a_i)=\a_i$.
Therefore $\th(\de-\a_i)\in \De\cup \{0\}$ and on the other hand
$\th(\de-\a_i)=-\de-\a_i<-\de$;
which causes a contradiction.

$(c)\Longrightarrow (d)$: Suppose $m_{\bar{\de}}>1$, then there exists $\a\in \De^+-\{\de\}$ such that $\bar{\a}=\bar{\de}$;
by the properties of root system, there exists $\be_1,\cdots,\be_k\in \Pi$ such that $\a=\de-\sum_{i=1}^k \be_i$ and
$$\de-\sum_{i=1}^j \be_i\in \De^+\qquad \mbox{for every }1\leq j\leq k.$$
Since $\bar{\a}=\bar{\de}$ and $\bar{\be}_i\in \Si^+\cup \{0\}$, we
have $\be_i\in \Pi_0$ for every $1\leq i\leq k$; especially $\be_1\in \Pi_0$ and $\be_1\notin \Pi\cap \de^\perp$; which contradict to $(c)$.

$(d)\Longrightarrow (b)$: Since
$$\ol{(-\de^\th)}=\f{-\de^\th+\th^2(\de)}{2}=\f{\de-\de^\th}{2}=\bar{\de}$$
and $m_{\bar{\de}}=1$, we have $\de^\th=-\de$ (by (1.8)). \textbf{Q.E.D}.
\newline

By Lemma 3.1, from (3.1) and those well known facts $|\sw{a}_l|=l(l+1)$, $|\sw{b}_l|=2l^2$, $|\sw{c}_l|=2l^2$,
$|\sw{d}_l|=2l(l-1)$, $|\sw{e}_6|=72$, $|\sw{e}_7|=126$, $|\sw{e}_8|=240$, $|\sw{f}_4|=48$, $|\sw{g}_2|=12$ (see \cite{Hel} p. 461-474),
we can obtain $(\de,\de)$ for every type of irreducible and reduced root systems:
\begin{eqnarray}
&&\De=\sw{a}_l:(\de,\de)=\f{1}{l+1};\ \De=\sw{b}_l:(\de,\de)=\f{1}{2l-1}; \ \De=\sw{c}_l:(\de,\de)=\f{1}{l+1};
\ \De=\sw{d}_l:(\de,\de)=\f{1}{2l-2};\nonumber\\
&&\De=\sw{e}_6:(\de,\de)=\f{1}{12};\ \De=\sw{e}_7:(\de,\de)=\f{1}{18};
\ \De=\sw{e}_8:(\de,\de)=\f{1}{30};\nonumber\\
&&\De=\sw{f}_4:(\de,\de)=\f{1}{9};\ \De=\sw{g}_2:(\de,\de)=\f{1}{4}.
\end{eqnarray}

Given $\sw{u},\Pi,\th,\Pi_0,\Si$, we can define the \textbf{Satake diagram} of $(\Pi,\th)$ as follows. Every root
of $\Pi_0$ is denoted by a black circle \begin{picture}(2,2)
\put(1,1){\circle*{2}}
\end{picture}\ \   and every root of $\Pi-\Pi_0$ by a white circle \begin{picture}(2,2)
\put(1,1){\circle{2}}
\end{picture}\ \  ;
if $\bar{\a}_i=\bar{\a}_j$ for $\a_i,\a_j\in \Pi-\Pi_0$, then $\a_i$ and $\a_j$ are joined by a curved
arrow. In \cite{Ar}, Araki gave the Satake diagram of $(\Pi,\th)$ and the Dynkin diagram of $\Si$
for all types of irreducible, simple and compact orthogonal involutive Lie algebras, i.e., $A\ I-A\ III,BD\ I,
C\ I-C\ II,D\ III,E\ I-E\ IX,F\ I-F\ II,G$. Then by Lemma 3.2, from the Satake diagram and (3.1), we can
justify whether $(\bar{\de},\bar{\de})=(\de,\de)$ or $(\bar{\de},\bar{\de})=1/2(\de,\de)$. For example, the Satake diagram
of $A\ II$ is
\setlength{\unitlength}{1mm}
\begin{center}
\begin{picture}(50,5)
\put(1.5,3.5){\circle*{2}}
\put(2.5,3.5){\line(1,0){5}}
\put(8.5,3.5){\circle{2}}
\put(9.5,3.5){\line(1,0){5}}
\put(15.5,3.5){\circle*{2}}
\put(16.5,3.5){\line(1,0){5}}
\put(21.8,2.6){$\cdots$}
\put(26,3.5){\line(1,0){5}}
\put(32,3.5){\circle{2}}
\put(33,3.5){\line(1,0){5}}
\put(39,3.5){\circle*{2}}
\put(0,0){$\a_1$}
\put(7,0){$\a_2$}
\put(14,0){$\a_3$}
\put(30,0){$\a_{l-1}$}
\put(38,0){$\a_l$}
\put(45,2.4){$l\geq 3$ is odd;}
\end{picture}
\end{center}
which yields $\Pi_0=\{\a_i:i\mbox{ is odd}\}$; thus $\a_1\in \Pi_0$ but $\a_1\notin \Pi\cap \de^\perp$ and Lemma 3.2
tells us $(\bar{\de},\bar{\de})=1/2(\de,\de)$. The ultimate results are: $(\bar{\de},\bar{\de})=1/2(\de,\de)$ when
$(\sw{u},\th)$ belongs to $A\ II,C\ II,E\ IV,F\ II$ or $(\sw{u},\th)$ belongs to $BD\ I$ and $l_1=1$;
otherwise $(\bar{\de},\bar{\de})=(\de,\de)$. Combining the results with (3.2), we can compute $(\bar{\de},\bar{\de})$,
i.e., $(\psi,\psi)$.
\newline

\textbf{Type II}. In this case, we denote $\sw{u}=\sw{v}\oplus \sw{v}$, where $\sw{v}$ is a compact and simple Lie
algebra; then $\th(X,Y)=(Y,X)$ for arbitrary $X,Y\in \sw{v}$, $\sw{k}_0=\{(X,X):X\in \sw{v}\}$, $\sw{p}_*
=\{(X,-X):X\in \sw{v}\}$. Let $\sw{t}$ be a maximal abelian subalgebra of $\sw{v}$, $\sw{t}_0=\sqrt{-1}\sw{t}$,
$\De^*\subset \sw{t}_0$ be the root system of $\sw{v}\otimes \C$ with respect to $\sw{t}\otimes \C$ with an ordering;
then $\sw{h}_{\sw{p}_*}
=\{(X,-X):X\in \sw{t}\}$ is a maximal abelian space of $\sw{p}_*$ and we can assume $\sw{h}_{\sw{k}_0}=\{(X,X):
X\in \sw{t}\}$; thus $\sw{h}_{\sw{p}_0}=\{(x,-x):x\in \sw{t}_0\}$, $\sw{h}_{\R}=\{(x,y):x,y\in \sw{t}_0\}$ and
\begin{eqnarray}
\De=(\De^*,0)\cup (0,\De^*),\qquad \Si=\{(\f{1}{2}\a,-\f{1}{2}\a):\a\in \De^*\}.
\end{eqnarray}
$\De$ has an lexicographic ordering induced by the ordering of $\De^*$, and we can define an ordering on $\Si$:
$(1/2\a,-1/2\a)>0$ if and only if $\a>0$; obviously $\De$ and $\Si$ have compatible orderings. Denote by $\de$
the highest root of $\De^*$, then $\psi=(1/2\de,-1/2\de)$ and
\begin{eqnarray}
(\psi,\psi)=\big((\f{1}{2}\de,-\f{1}{2}\de),(\f{1}{2}\de,-\f{1}{2}\de)\big)=\f{1}{2}(\de,\de);
\end{eqnarray}
i.e., the square of the length of the highest restricted root is a half of the square of the length of the highest root of $\De^*$.

\section{Computation of injectivity radius and diameter}
\setcounter{equation}{0}

In this section, we assume $U=G/K$ be irreducible; it is well known there is a one-to-one correspondence between all
compact, simply connected and irreducible Riemannian symmetric spaces and all
semisimple, compact, reduced and irreducible orthogonal involutive Lie algebras (cf. \cite{Hel} p. 438-443).
Denote by $(\sw{u},\th,\lan,\ran)$ the orthogonal involutive Lie algebras corresponding to $U=G/K$, then
\begin{eqnarray}
\lan,\ran=-\ep(,)
\end{eqnarray}
for some positive instant $\ep$ (cf. \cite{AB} p. 23-26).

\begin{Rem}
$\ep$ has  geometrical significance: there is a close relationship between $\ep$ and the Ricci curvature of $M$.
\end{Rem}

The denotation of $\sw{k}_0,\sw{p}_*$ is similar to Section 1. Let $\{Y_1,\cdots,Y_m\}$ be a basis of $\sw{k}_0$,
$\{X_1,\cdots,X_n\}$ be a basis of $\sw{p}_*$; (1.4) tells us $[\sw{k}_0,\sw{p}_*]\subset \sw{p}_*$ and
$[\sw{p}_*,\sw{p}_*]\subset \sw{k}_0$, then for arbitrary
$X\in \sw{p}_*$,  the matrix of $ad X|_{\sw{u}}$ with respect to the basis $\{Y_1,\cdots,Y_m,X_1,\cdots,X_n\}$ is
$$\left(\begin{array}{cc} & C\\D & \end{array}\right);\qquad (C\mbox{ is a }m\times n\mbox{ matrix},D\mbox{ is a }n\times m\mbox{ matrix})$$
then
$$(ad X)^2=\left(\begin{array}{cc} CD & \\ & DC\end{array}\right)$$
and therefore
\begin{eqnarray}
(X,X)=\tr(ad X)^2|_{\sw{u}}=\tr(CD)+\tr(DC)=2\ \tr(DC)=2\ \tr(ad X)^2|_{\sw{p}_*}.
\end{eqnarray}
By (1.9), we have
\begin{eqnarray}
Ric(X,X)&=&R_{X,X_i}X,X_i=\lan ad[X,X_i]X,X_i\ran=\lan -(ad\ X)^2 X_i,X_i\ran\nonumber\\
&=&-\tr (ad\ X)^2|_{\sw{p}_*}=-\f{1}{2}(X,X)=\f{1}{2\epsilon}\lan X,X\ran;
\end{eqnarray}
i.e., $M$ is an Einstein manifold with Ricci curvature $1/(2\epsilon)$.
\newline

By the definition of injectivity radius and diameter, $i(M)=\min_{p\in M,X\in C(p)}|X|$ and $d(M)=\max_{p\in M,X\in C(p)}
|X|$; by Theorem 3.1, for every $X\in C(p)$ where $p=aK$, we have $X=(dL_a)_o \big(Ad(k)(\pi\sqrt{-1}x\big)\big)$
for some $k\in K$ and $x\in \triangle'$; which yields $|X|=\pi\ep^{1/2}(x,x)^{1/2}$ and
\begin{eqnarray}
i(M)=\pi\ep^{1/2}i(\Si),\qquad d(M)=\pi\ep^{1/2}d(\Si).
\end{eqnarray}
($i(\Si)$ and $d(\Si)$ is defined in (2.2), notice that the inner product $(,)$ on $\Si$ is induced by the Killing form
of $\sw{u}$.)

According to (4.4), we have the following theorem on $i(M)$:

\begin{Thm}
Let $M$ be a simply-connected, compact and irreducible Riemannian symmetric space, $\k$ be the maximum of the sectional curvatures of
$M$, then $i(M)=\pi\k^{-1/2}$.

\end{Thm}

\textbf{Proof}: When $\lan,\ran=-(,)$, i.e., $\ep=1$, it is known that $\k=(\psi,\psi)$, where $\psi$ denotes the highest
restricted root (see \cite{Hel} p.334). Then for general cases such that $\ep\neq 1$, we have $\k=\ep^{-1}(\psi,\psi)$.
On the other hand, by (4.4) and (2.3), $i(M)=\pi\ep^{1/2}(\psi,\psi)^{-1/2}$. Thus $i(M)=\pi\k^{-1/2}$.
\textbf{Q.E.D}.
\newline

Moreover, from the results obtained in Section 2 and Section 3, we can compute $i(M)$ and $d(M)$ for every type
of compact, simply connected and irreducible Riemannian symmetric spaces and list the results in Table 4.1 and Table 4.2.

\newcommand{\ZZ}[2]{\rule[#1]{0pt}{#2}}
\newcommand{\rb}[1]{\raisebox{#1}[0pt]}
\begin{center}
\begin{tabular}{|c|c|c|c|c|c|}
\multicolumn{6}{c}{\textbf{Table 4.1}}\\
\multicolumn{6}{c}{}\\
\multicolumn{6}{c}{\textit{The injectivity radius and diameter of compact, simply connected and irreducible}}\\
\multicolumn{6}{c}{ \textit{Riemannian symmetric spaces of Type I when $\ep=1$, i.e., $Ric=1/2$}}\\
\multicolumn{6}{c}{}\\
\hline
\ZZ{-8pt}{22pt}\textbf{ Type} & $ M$ & $\Si$   & $(\bar{\de},\bar{\de})$ & $i(M)$  & $d(M)$ \\\hline
\rb{-1.8ex}{$A\ I$}  &  \rb{-1.8ex}{$SU(n)/SO(n)$}  & \rb{-1.8ex}{$\sw{a}_{n-1}$}    &  \rb{-1.8ex}{$\f{1}{n}$}   &  \rb{-1.8ex}{$\pi n^{1/2}$}  & \ZZ{-1pt}{15pt}$\f{\sqrt{2}}{2}\pi n$\ ($n$ is even)\\
                     &                              &                                &                            &                                       & \rb{1ex}{$\f{\sqrt{2}}{2}\pi(n^2-1)^{1/2}$\ ($n$ is odd)}\\ \hline
\rb{-1.8ex}{$A\ II$} &  \rb{-1.8ex}{$SU(2n)/Sp(n)$} & \rb{-1.8ex}{$\sw{a}_{n-1}$}    &  \rb{-1.8ex}{$\f{1}{4n}$}  &  \rb{-1.8ex}{$2\pi n^{1/2}$} & \ZZ{-1pt}{12pt}$\sqrt{2}\pi n$ ($n$ is even)\\
                     &                              &                                &                                         &                                       & \rb{1ex}{$\sqrt{2}\pi(n^2-1)^{1/2}$\ ($n$ is odd)}\\ \hline
                     &                              & $(\sw{bc})_p$     &                                         &                                  &                                \\
                     &                              & ($2\leq p<q$)     &                                         &                                  &                                 \\
$A\ III$             &  $G_{p,q}(\C)$                & $\sw{c}_p$        &  $\f{1}{p+q}$              &   $\pi(p+q)^{1/2}$      &  $\pi(p+q)^{1/2}p^{1/2}$    \\
                     &  ($p\leq q$)                 & ($p=q\geq 2$)     &                                         &                                  &                          \\
                     &                              & $(\sw{bc})_1$     &                                         &                                  &                           \\
                     &                              &  $(p=1)$          &                                         &                                  &                           \\\hline
\ZZ{-10pt}{26pt}$C\ I$               &  $SP(n)/U(n)$                & $\sw{c}_n$        &  $\f{1}{n+1}$       &  $\pi(n+1)^{1/2}$                  & $\pi(n+1)^{1/2}n^{1/2}$\\ \hline
                     &                              & $(\sw{bc})_p$     &                                         &                                  &                                \\
                     &                              & ($2\leq p<q$)     &                                         &                                  &                                 \\
$C\ II$              &  $G_{p,q}(\H)$                & $\sw{c}_p$        &  $\f{1}{2(p+q+1)}$         &   $\sqrt{2}\pi(p+q+1)^{1/2}$      &  $\sqrt{2}\pi(p+q+1)^{1/2}p^{1/2}$    \\
                     &  ($p\leq q$)                 & ($p=q\geq 2$)     &                                         &                                  &                          \\
                     &                              & $(\sw{bc})_1$     &                                         &                                  &                           \\
                     &                              &  $(p=1)$          &                                         &                                  &                           \\\hline
                     &                              &                   &                                         &                                      & $\sqrt{2}\pi(p+q-2)^{1/2}$   \\
                     &                              & $\sw{b}_p$        &  $\f{1}{p+q-2}$              &   $\pi(p+q-2)^{1/2}$        & ($p\leq 3$)\\
                     &                              & ($2\leq p<q$)     &                                         &                                      & $\f{\sqrt{2}}{2}\pi(p+q-2)^{1/2}p^{1/2}$\\
$BD\ I$              &  $G_{p,q}(\R)$                &                   &                                         &                                      & ($p\geq 4$)\\\cline{3-6}
                     &  ($p\leq q$)                 & $\sw{d}_p$        &  \rb{-1.8ex}{$\f{1}{2p-2}$}  &   \rb{-1.8ex}{$\sqrt{2}\pi(p-1)^{1/2}$}  & \rb{-1.8ex}{$\pi(p-1)^{1/2}p^{1/2}$} \\
                     &                              & \rb{1.8ex}{($4\leq p=q$)}     &                                         &                                      &\\\cline{3-6}
                     &                              & $\sw{a}_1$        &  \rb{-1.8ex}{$\f{1}{2q-2}$}    &   \rb{-1.8ex}{$\sqrt{2}\pi(q-1)^{1/2}$}  & \rb{-1.8ex}{$\sqrt{2}\pi(q-1)^{1/2}$}\\
                     &                              & \rb{1.8ex}{($1=p<q$)}         &                                         &                                      & \\\hline
\end{tabular}
\end{center}

\begin{center}
\begin{tabular}{|c|c|c|c|c|c|}
\multicolumn{6}{c}{\textbf{Table 4.1}(\textit{continued})}\\
\multicolumn{6}{c}{}\\
\hline
\ZZ{-8pt}{22pt}\textbf{ Type} & $ M$ & $\Si$   & $(\bar{\de},\bar{\de})$ & $i(M)$  & $d(M)$ \\\hline
              &                            &    $\sw{c}_{\f{n}{2}}$  &  \rb{-1.8ex}{$\f{1}{2n-2}$}        &   \rb{-1.8ex}{$\sqrt{2}\pi(n-1)^{1/2}$}   & \rb{-1.8ex}{$\pi(n-1)^{1/2}n^{1/2}$}\\
$D\ III$       &   $SO(2n)/U(n)$            &    \rb{1.2ex}{($n$ is even)}           &                                    &                                       &                               \\\cline{3-6}
               &                            &    $(\sw{bc})_{\f{n-1}{2}}$ & \rb{-1.8ex}{$\f{1}{2n-2}$}       &   \rb{-1.8ex}{$\sqrt{2}\pi(n-1)^{1/2}$}   & \rb{-1.8ex}{$\pi(n-1)$}\\
               &                            &    \rb{1.2ex}{($n$ is odd)}              &                                  &                                       &                    \\\hline
\ZZ{-9pt}{24pt}$E\ I$         &   $(\sw{e}_6,sp(4))$       &    $\sw{e}_6$              &  $\f{1}{12}$       &  $2\sqrt{3}\pi$       &      $4\sqrt{2}\pi$\\\hline
\ZZ{-9pt}{24pt}$E\ II$        &   $(\sw{e}_6,su(6)\oplus su(2))$ & $\sw{f}_4$  &   $\f{1}{12}$      &  $2\sqrt{3}\pi$     & $2\sqrt{6}\pi$\\\hline
\ZZ{-9pt}{24pt}$E\ III$       &   $(\sw{e}_6,so(10)\oplus \R)$    & $(\sw{bc})_2$ & $\f{1}{12}$     & $2\sqrt{3}\pi$   &   $2\sqrt{6}\pi$\\\hline
\ZZ{-9pt}{24pt}$E\ IV$        &   $(\sw{e}_6,\sw{f}_4)$     &   $\sw{a}_2$     & $\f{1}{24}$      & $2\sqrt{6}\pi$    &   $4\sqrt{2}\pi$\\\hline
\ZZ{-9pt}{24pt}$E\ V$         &   $(\sw{e}_7,su(8))$        &   $\sw{e}_7$     & $\f{1}{18}$  &  $3\sqrt{2}\pi$       &  $3\sqrt{6}\pi$\\\hline
\ZZ{-9pt}{24pt}$E\ VI$         &   $(\sw{e}_7,so(12)\oplus su(2))$        &   $\sw{f}_4$     & $\f{1}{18}$  &  $3\sqrt{2}\pi$       &  $6\pi$\\\hline
\ZZ{-9pt}{24pt}$E\ VII$    &   $(\sw{e}_7,\sw{e}_6\oplus \R)$   &  $\sw{c}_3$  &  $\f{1}{18}$     &  $3\sqrt{2}\pi$     &  $3\sqrt{6}\pi$\\\hline
\ZZ{-9pt}{24pt}$E\ VIII$         &   $(\sw{e}_8,so(16))$        &   $\sw{e}_8$     & $\f{1}{30}$  &  $\sqrt{30}\pi$       &  $2\sqrt{15}\pi$\\\hline
\ZZ{-9pt}{24pt}$E\ IX$         &   $(\sw{e}_8,\sw{f}_7\oplus su(2))$        &   $\sw{f}_4$     & $\f{1}{30}$  &  $\sqrt{30}\pi$       &  $2\sqrt{15}\pi$\\\hline
\ZZ{-9pt}{24pt}$F\ I$         &   $(\sw{f}_4,sp(3)\oplus su(2))$        &   $\sw{f}_4$     & $\f{1}{9}$  &  $3\pi$       &  $3\sqrt{2}\pi$\\\hline
\ZZ{-9pt}{24pt}$F\ II$         &   $(\sw{f}_4,so(9))$        &   $(\sw{bc})_1$     & $\f{1}{18}$  &  $3\sqrt{2}\pi$       &  $3\sqrt{2}\pi$\\\hline
\ZZ{-9pt}{24pt}$G$         &   $(\sw{g}_2,su(2)\oplus su(2))$        &   $\sw{g}_2$     & $\f{1}{4}$  &  $2\pi$       &  $\f{4\sqrt{3}}{3}\pi$\\\hline
\end{tabular}
\end{center}

\begin{center}
\begin{tabular}{|c|c|c|c|c|}
\multicolumn{5}{c}{\textbf{Table 4.2}}\\
\multicolumn{5}{c}{}\\
\multicolumn{5}{c}{\textit{The injectivity radius and diameter of compact, simply connected and irreducible}}\\
\multicolumn{5}{c}{ \textit{Riemannian symmetric spaces of Type II when $\ep=1$, i.e., $Ric=1/2$}}\\
\multicolumn{5}{c}{}\\
\hline
\ZZ{-8pt}{22pt} $M$ & $\De^*$   & $(\de,\de)$ & $i(M)$  & $d(M)$ \\\hline
\rb{-1.8ex}{$SU(n)$}     &  \rb{-1.8ex}{$\sw{a}_{n-1}$}   & \rb{-1.8ex}{$\f{1}{n}$}  & \rb{-1.8ex}{$\sqrt{2}\pi n^{1/2}$} & $\pi n$($n$ is even)\\
            &                   &             &                       & \rb{1.2ex}{$\pi(n^2-1)^{1/2}$($n$ is odd)}\\\hline
\rb{-1.8ex}{$Spin(2n+1)$}&  \rb{-1.8ex}{$\sw{b}_n$}       & \rb{-1.8ex}{$\f{1}{2n-1}$} & \rb{-1.8ex}{$\sqrt{2}\pi(2n-1)^{1/2}$} & $2\pi(2n-1)^{1/2}$($n\leq 3$)\\
            &                   &               &                           & \rb{1.2ex}{$\pi(2n-1)^{1/2}n^{1/2}$($n\geq 4$)}\\\hline
\ZZ{-9pt}{24pt}$Sp(n)$     & $\sw{c}_n$       & $\f{1}{n+1}$   & $\sqrt{2}\pi(n+1)^{1/2}$  & $\sqrt{2}\pi(n+1)^{1/2}n^{1/2}$\\\hline
\ZZ{-9pt}{24pt}$Spin(2n)$  & $\sw{d}_n$       & $\f{1}{2n-2}$  & $2\pi(n-1)^{1/2}$         & $\sqrt{2}\pi(n-1)^{1/2}n^{1/2}$\\\hline
\ZZ{-9pt}{24pt}$E_6$       & $\sw{e}_6$       & $\f{1}{12}$    & $2\sqrt{6}\pi$       &   $8\pi$\\\hline
\ZZ{-9pt}{24pt}$E_7$       & $\sw{e}_7$       & $\f{1}{18}$    & $6\pi$               &   $6\sqrt{3}\pi$\\\hline
\ZZ{-9pt}{24pt}$E_8$       & $\sw{e}_8$       & $\f{1}{30}$    & $2\sqrt{15}\pi$      &   $2\sqrt{30}\pi$\\\hline
\ZZ{-9pt}{24pt}$F_4$       & $\sw{f}_4$       & $\f{1}{9}$     & $3\sqrt{2}\pi$       &   $6\pi$\\\hline
\ZZ{-9pt}{24pt}$G_2$       & $\sw{g}_2$       & $\f{1}{4}$     & $2\sqrt{2}\pi$       &   $\f{4\sqrt{6}}{3}\pi$\\\hline
\end{tabular}
\end{center}

\begin{Rem}
In Table 4.1, $\Si$ denotes the restricted root system, $\bar{\de}$ denotes the highest restricted root, $(\bar{\de},\bar{\de})$
denotes the square of the length of $\bar{\de}$, $i(M)$ and $d(M)$ denote injectivity radius and diameter of $M$,
respectively. In Table 4.2, $M$ is a compact, simply connected and simple Lie group with bi-invariant metric;
$\sw{v}$ is the Lie algebra associated to $M$, $\sw{t}$ is a maximal abelian subalgebra of $\sw{v}$
and $\De^*$ denotes the root system of $\sw{v}\otimes \C$ with respect to $\sw{t}\otimes \C$, $\de$ denotes the highest root
of $\De^*$ (cf. Section 3).
\end{Rem}

\begin{Rem}
In Table 4.1 and Table 4.2, we assume $\ep=1$, i.e., the $K$-invariant metric on $M=U/K$ is induced by $-(,)$ on
$\sw{u}$, and $Ric=1/2$. For general cases such that $\ep\neq 1$, we should multiply the corresponding results in Table 4.1
or Table 4.2 by $\ep^{1/2}$.
\end{Rem}

For example, let $U=SO(p+q)$, $K=SO(p)\times SO(q)$, $M=U/K=G_{p,q}(\R)$ $(p\leq q)$; then $\sw{u}=so(p+q)$, $\sw{k}_0=so(p)\oplus so(q)$ and
\begin{eqnarray}
\sw{p}_*=\Big\{\left(\begin{array}{cc}0 & -X^T\\X & 0\end{array}\right)\in so(p+q):X\mbox{ is a }q\times p\mbox{ matrix}\Big\}.
\end{eqnarray}
Denote $\lan A,B\ran=-1/2\tr(AB)$ for every $A,B\in \sw{p}_*$, then it is easily seen that $\lan,\ran$ is invariant under
$K$; the canonical metric on $G_{p,q}(\R)$ is induced by $\lan,\ran$ and it's $U$-invariant (cf. \cite{Ko2} p.271-273). It is well
known that $(A,B)=(p+q-2)\tr(AB)$ for every $A,B\in so(p+q)$; so $\ep=1/\big(2(p+q-2)\big)$ and furthermore, from
Table 4.1 we have
\begin{eqnarray}
i\big(G_{p,q}(\R)\big)=\left\{\begin{array}{cc} \f{\sqrt{2}}{2}\pi & p\geq 2;\\ \pi & p=1.\end{array}\right.
\qquad d\big(G_{p,q}(\R)\big)=\left\{\begin{array}{cc} \pi & p=1\mbox{ or }2\leq p\leq 3\mbox{ and }q>p;\\\f{1}{2}\pi p^{1/2} & \mbox{ otherwise}.\end{array}\right.
\end{eqnarray}

\begin{Rem}
If $(M,g)$ is a compact, simply connected and reducible Riemannian symmetric space, then by de Rham decomposition Theorem
(see \cite{Ko2} p. 210-216), $(M,g)=(M_1,g_1)\times\cdots\times (M_r,g_r)$, where $(M_1,g_1),\cdots,(M_r,g_r)$ are all compact, simply connected and
irreducible. For every $p=(p_1,\cdots,p_r)\in M$ and $X=(X_1,\cdots,X_r)\in T_p M$, where $p_i\in M_i$, $X_i\in T_{p_i}M_i$,
$d(p,\exp_p(X))=|X|$ if and only if $d(p_i,\exp_{p_i}(X_i))=|X_i|$ for every $1\leq i\leq r$;
therefore
\begin{eqnarray}
i(M)=\min_{1\leq i\leq r}i(M_i),\qquad d(M)=\big(\sum_{i=1}^r d(M_i)^2\big)^{1/2}.
\end{eqnarray}
\end{Rem}


\begin{thebibliography}{AB}
\bibitem[Ar]{Ar}Araki, S. I, \textit{On root systems and an infinitesimal classification of irreducible symmetric spaces}, J.Math. Osaka City Univ.\textbf{13}(1962), 1-34.
\bibitem[Bo]{AB}A. Borel, \textit{Semisimple Groups and Riemannian Symmetric Spaces}, Hindustan Book Agency(1998).
\bibitem[CE]{CE}J. Cheeger and D.Ebin, \textit{Comparison Theorem in Riemannian Geometry}, Noth-Holland publishing company(1975).
\bibitem[Ch]{Ch}J. Cheeger, \textit{Pinching theorems for a certain class of Riemannian manifolds}, Amer. J. Math. \textbf{91}(1969),807-834.
\bibitem[Cr]{Cr}R. Crittenden, \textit{Minimum and conjugate points in symmetric spaces}, Canadian J. Math. \textbf{14}(1962),320-328.
\bibitem[GS]{GS}Grove, Karsten and Shiohama, Katsuhiro
\textit{A generalized sphere theorem,}
Ann. Math. (2) \textbf{106} (1977), no. 2, 201--211.
\bibitem[Hel]{Hel} S. Helgason, \textit{Differential Geometry, Lie groups, and Symmetric Spaces}, Academic Press(1978).
\bibitem[Kna]{Kna} A.W.Knapp, \textit{Lie group, Beyond an introduction}, Progress in mathematics(2002).
\bibitem[Ko1]{Ko}S. Kobayashi and K. Nomizu, \textit{Foundations of differential geometry}, Volume 1, Interscience Publishers(1963).
\bibitem[Ko2]{Ko2}S. Kobayashi and K. Nomizu, \textit{Foundations of differential geometry}, Volume 2, Interscience Publishers(1963).
\bibitem[Xin]{Xin}Yuanlong Xin, \textit{Minimal Submanifolds and Related Topics}, Nankai Tracts in Mathematics(2003).
\end{thebibliography}
\end{document}